\documentclass[twoside,reqno]{amsart} 

\usepackage{array}
\usepackage{times}
\usepackage{amsthm}           
\usepackage{enumerate}        
\usepackage{euscript}         
\usepackage{dan2e}

\numberwithin{equation}{section}

\newenvironment{roenumerate}{\begin{enumerate}[\ulp\itshape i\urp]}{\end{enumerate}}

\newtheorem{thm}[equation]{Theorem} 
\newtheorem{defn}[equation]{Definition}
\newtheorem{prop}[equation]{Proposition}
\newtheorem{cor}[equation]{Corollary}
\newtheorem{lemma}[equation]{Lemma}

\theoremstyle{definition}  
\newtheorem{example}[equation]{Example}

\newcommand{\ulp}{\textup{(}}
\newcommand{\urp}{\textup{)}}

\newcommand{\dfn}{\textbf} 
\newcommand{\mdfn}[1]{\dfn{\mathversion{bold}#1}} 

\newcommand{\cat}               {\EuScript}  
\renewcommand{\cA}              {{\cat A}}
\newcommand{\sA}              {{s\cA}}
\renewcommand{\cB}              {{\cat B}}
\renewcommand{\cC}              {{\cat C}}
\renewcommand{\cD}              {{\cat D}}

\renewcommand{\cM}              {{\cat M}}
\renewcommand{\cP}              {{\cat P}}

\renewcommand{\cS}              {{\cat S}}

\newcommand{\A}                 {{\cat A}}
\newcommand{\Ch}                {\textup{Ch}}

\newcommand{\D}                 {\cD}

\newcommand{\DC}                {\D_{\cC}}
\newcommand{\DP}                {\D_{\cP}}

\DeclareMathOperator{\Ex}     {Ex}

\newcommand{\Proj}              {\cP}

\newcommand{\Exact}             {\cE}

\DeclareMathOperator{\Pure}     {P}
\renewcommand{\PExt}              {\Pure\!\Ext}

\newcommand{\Iinj}              {\textup{$I$-inj}}
\newcommand{\Jinj}              {\textup{$J$-inj}}
\newcommand{\Icof}              {\textup{$I$-cof}}
\newcommand{\Jcof}              {\textup{$J$-cof}}
\newcommand{\Icell}              {\textup{$I$-cell}}
\newcommand{\Jcell}              {\textup{$J$-cell}}

\renewcommand{\st}{\,\mid\,}

\begin{document}

\title{Derived categories and projective classes}
\author{J.\ Daniel Christensen}

\address{
Department of Mathematics\\
Johns Hopkins University\\
3400 N.~Charles St.\\
Baltimore, MD 21218-2686}
\email{jdc@math.jhu.edu}



\begin{abstract}
An important example of a model category is the category of unbounded
chain complexes of $R$-modules,
which has as its homotopy category the derived category of the ring $R$.
This example shows that traditional homological algebra is encompassed
by Quillen's homotopical algebra.
The goal of this paper is to show that more general forms of homological 
algebra also fit into Quillen's framework.
Specifically, any set of objects in a complete and cocomplete
abelian category $\cA$ generates a projective class on
$\cA$, which is exactly the information needed to do homological
algebra in $\cA$.  
The main result is that if the generating objects are ``small''
in an appropriate sense,
then the category of chain complexes of objects of $\cA$ has a 
model category structure which reflects the homological algebra
of the projective class.
The motivation for this work is the construction of the ``pure
derived category'' of a ring $R$.
Finally, we explain how the category of simplicial objects in 
a possibly non-abelian category can be equipped with a model category
structure reflecting a given projective class.
\end{abstract}

\date{May 9, 2000}

\maketitle

\tableofcontents

\renewcommand{\baselinestretch}{1.15}\normalsize

\section*{Introduction}\label{se:intro}

An important example of a model category is the category $\Ch$ of unbounded
chain complexes of $R$-modules,
which has as its homotopy category the derived category $\D$ of the ring $R$.
The formation of a projective resolution is an example of cofibrant
replacement, and traditional derived functors are examples of
derived functors in the model category sense.
This example shows that traditional homological algebra is encompassed
by Quillen's homotopical algebra, and indeed this unification was one
of the main points of Quillen's work~\cite{qu:ha}.

The goal of this paper is to illustrate that more general forms of 
homological algebra also fit into Quillen's framework.
In any abelian category $\A$ there is a natural notion of ``projective
object'' and ``exact sequence.''  
However, it is sometimes useful to impose different definitions
of these terms.
If this is done in a way that satisfies some natural axioms,
what is obtained is a ``projective class,'' which is exactly 
the information needed to do homological algebra in $\A$.  
Our main result is that given any projective class satisfying 
a set-theoretic hypothesis,
the category of chain complexes of objects of $\A$ has a 
model category structure which reflects the homological algebra
of the projective class.
The motivation for this work is the construction of the ``pure
derived category'' of a ring $R$.
Pure homological algebra has applications to phantom maps in
the stable homotopy category and in the (usual) derived
category of a ring, and these connections will be described.

When $\A$ has enough projectives, the projective objects and
exact sequences form a projective class.  Therefore the results
of this paper apply to traditional homological algebra as well.
Even in this special case, it is not a trivial fact that the
category of \emph{unbounded} chain complexes can be given a model 
category structure, and indeed Quillen restricted himself to the
bounded below case.
I know of three other written proofs that the category of unbounded
chain complexes is a model category~\cite{gr:mc,hi:haha,ho:mc}, which do the
case of $R$-modules, but this was probably known to others as well.

An important corollary of the fact that a derived category $\D$ is the
homotopy category of a model category is that $\D(X,Y)$ is a set (as
opposed to a proper class) for any two chain complexes $X$ and $Y$.
This is not the case in general, and much work on derived categories 
ignores this possibility.  
It is not just a pedantic point;  if one uses the morphisms in the
derived category to index constructions in other categories or to
define cohomology groups, one needs to know that the indexing class
is actually a set.
Recently, the unbounded case has been handled 
under various assumptions on $\A$.
(See Weibel~\cite{we:iha} Remark 10.4.5, which credits Gabber, and
Exercise 10.4.5, which credits Lewis, May and Steinberger~\cite{lemast:esht}.
See also Kriz and May~\cite[Part III]{krma:oamm}.)
The assumptions that appear in the present paper are different from those that
have appeared before and the proof is somewhat easier (because
of our use of the theory of cofibrantly generated model categories),
so this paper may be of some interest even in this special case.

Another consequence of the fact that $\Ch$ is a model category is
the existence of resolutions coming from cofibrant and fibrant
approximations, and the related derived functors.  
Some of these are discussed in \cite{avfoha} and \cite{sp:ruc}.
We do not discuss these topics here, but just mention that
the model category point of view can provide a framework for
some of this material.

We also briefly discuss the category of non-negatively graded
chain complexes.  In this case we describe a model category structure
that works for an arbitrary projective class, without any set-theoretic 
hypotheses.
More generally, we show that under appropriate hypotheses a projective
class on a possibly non-abelian category $\cA$ determines a model
category structure on the category of simplicial objects in $\cA$.
As an example, we deduce that the category of equivariant simplicial
sets has various model category structures.

We now briefly outline the paper.  In Section~\ref{se:pc} we
give the axioms for a projective class and mention a few examples.
In Section~\ref{se:dc} we describe the model category structure
on the category of chain complexes which takes into account a
given projective class.  Then we state the set-theoretic assumption
needed and prove our main result, using the recognition lemma for
cofibrantly generated categories, which is recalled in Section~\ref{se:mc}.
In Section~\ref{se:examples} we give two examples, the traditional
derived category of $R$-modules and the pure derived category,
and we describe why the pure derived category is interesting.
In the final section we discuss the bounded below case,
which works for any projective class, and describe a result for
simplicial objects in a possibly non-abelian category.

I thank Haynes Miller for asking the question which led to this paper
and Mark Hovey, Haynes Miller and John Palmieri for fruitful and
enjoyable discussions.

\section{Projective classes}\label{se:pc}

In this section we explain the notion of a projective class, which is
the information necessary in order to do homological algebra.
Intuitively, a projective class is a choice of which sort of 
``elements'' we wish to think about.

The elements of a set $X$ correspond bijectively to the
maps from a singleton to $X$, and
the elements of an abelian group $A$ correspond bijectively to
the maps from $\Z$ to $A$.
Motivated by this, we 
call a map $P \ra A$ in any category a \mdfn{$P$-element of $A$}.
If we don't wish to mention $P$, we call such a map a \mdfn{generalized
element of $A$}.
A map $A \ra B$ in any category is determined by what it does on
generalized elements.
If $\Proj$ is a collection of objects, then a \mdfn{$\Proj$-element} means
a $P$-element for some $P$ in $\Proj$.

Let $\A$ be a \dfn{pointed category}, \ie a category in which initial and
terminal objects exist and agree.  In a pointed category, any
initial (equivalently, terminal) object is called a \dfn{zero object}.
Let $P$ be an object of $\A$.  
A sequence 
\[
A \lra B \lra C
\]
is said to be \mdfn{$P$-exact} if the composite $A \ra C$ is the zero map 
(the unique map which factors through a zero object) and
\[
\A(P,A) \lra \A(P,B) \lra \A(P,C)
\]
is an exact sequence of pointed sets (the base points being the zero maps).
The latter can be rephrased as the condition that $A \ra B \ra C$
induces an exact sequence of $P$-elements.
A \mdfn{$\Proj$-exact sequence} is one which is $P$-exact for all
$P$ in $\Proj$.
A map $A \ra B$ is \mdfn{$P$-epic} (resp.\ \mdfn{$\Proj$-epic}) if it induces a
surjection of $P$-elements (resp.\ $\Proj$-elements).

A \dfn{projective class} on $\A$ is a collection $\Proj$ of objects
of $\A$ and a collection $\Exact$ of sequences $A \ra B \ra C$ in $\A$ 
such that 
\begin{roenumerate}
\item $\Exact$ is precisely the collection of all $\Proj$-exact sequences;
\item $\Proj$ is precisely the collection of all objects $P$ such that each
      sequence in $\Exact$ is $P$-exact;
\item any map $A \ra B$ can be extended to a sequence $P \ra A \ra B$
      in $\Exact$ with $P$ in $\Proj$.
\end{roenumerate}
When a collection $\Proj$ is part of a projective class $(\Proj,\Exact)$,
the projective class is unique, and so we say that $\Proj$ determines
a projective class or even that $\Proj$ is a projective class.
An object of $\Proj$ is called a \mdfn{$\Proj$-projective}.
Condition (iii) says that there are ``enough $\Proj$-projectives.''

\begin{example}
Let $\A$ be the category of pointed sets, let $\Proj$ be the collection
of all pointed sets and let $\Exact$ be the collection of all exact sequences 
of pointed sets.  Then $\Exact$ is precisely the collection of $\Proj$-exact
sequences, and $\Proj$ is a projective class.
\end{example}

\begin{example}\label{ex:cp}
For an associative ring $R$,
let $\A$ be the category of left $R$-modules, let $\Proj$ be the collection
of all summands of free $R$-modules and let $\Exact$ be the collection of
all exact sequences of $R$-modules.  Then $\Exact$ is precisely the
collection of $\Proj$-exact sequences, and $\Proj$ is a projective
class.
\end{example}

The above two examples are ``categorical'' projective classes in
the sense that the $\Proj$-epi\-mor\-phisms are just the epimorphisms
and the $\Proj$-projectives are the categorical projectives,
\ie those objects $P$ such that maps from $P$ lift through epimorphisms.

Here are two examples of non-categorical projective classes.

\begin{example}\label{ex:pp}
Let $\A$ be the category of left $R$-modules, as in Example~\ref{ex:cp}.
Let $\Proj$ consist of all summands of sums of finitely presented
modules and define $\Exact$ to consist of all $\Proj$-exact sequences.
Then $\Proj$ is a projective class.
A short exact sequence is $\Proj$-exact iff it remains exact after
tensoring with every right module.
\end{example}

\begin{example}\label{ex:ph}
Let $\cS$ be the homotopy category of spectra and let $\Proj$ 
consist of all retracts of wedges of finite spectra.
Then $\Proj$ determines a projective class, and a map is a $\Proj$-epimorphism
iff its cofibre is a phantom map.
\end{example}

Examples~\ref{ex:cp} and~\ref{ex:pp} will be discussed further in
Section~\ref{se:examples}.
Example~\ref{ex:ph} is studied in~\cite{ch:itcpgs}, along with
similar examples.

We note that all of the examples mentioned above are \dfn{determined
by a set} in the sense of the following lemma, whose proof is easy.

\begin{lemma}\label{le}
Suppose $\cF$ is any \emph{set} of objects in a pointed category with
coproducts.
Let $\Exact$ be the collection of $\cF$-exact sequences and 
let $\Proj$ be the collection of all objects $P$ such that every
sequence in $\Exact$ is $P$-exact.
Then $\Proj$ is the collection of retracts of coproducts of objects
of $\cF$ and $(\Proj,\Exact)$ is a projective class. \qed
\end{lemma}

A projective class is precisely the information needed to form
projective resolutions and define derived functors.
For further information, we refer the reader to~\cite{eimo:frha}
for the general theory and to~\cite{ch:itcpgs} in the special
case of a triangulated category.

\section{Cofibrantly generated model categories}\label{se:mc}

In this section we briefly recall the basics of cofibrantly generated
model categories.
This material will be used in the next section to prove our main result.
For more details, see the books by Dwyer, Hirschhorn and Kan~\cite{dwhika}
and Hirschhorn~\cite{hi:lmc}.
We assume knowledge of the basics of model categories, for 
which~\cite{dwsp:htmc} is an excellent reference.

We will always assume our model categories to be complete and cocomplete.

In the following, a cardinal number is thought of as the first
ordinal with that cardinality.

\begin{defn}
  Given an ordinal $\kappa$, a \mdfn{$\kappa$-sequence in a category $\cM$}
  is a diagram  
\[ 
X_{0} \lra X_{1} \lra \cdots \lra X_{\lambda} \lra X_{\lambda + 1} \lra \cdots 
\]
  indexed by the ordinals less than $\kappa$, such that for each limit
  ordinal $\gamma$ less than $\kappa$ the natural map 
  $\colim_{\lambda < \gamma} X_{\lambda} \ra X_{\gamma}$
  is an isomorphism.

  A cardinal $\gamma$ is said to be \dfn{regular} if any sum of fewer
  than $\gamma$ cardinals each less than $\gamma$ is less than $\gamma$.
  The first few regular cardinals are $0$, $1$, $2$ and $\aleph_{0}$.

  An object $P$ is said to be \mdfn{small relative to a subcategory $\cB$} 
  if there exists a cardinal $\kappa$ such that for any regular
  cardinal $\gamma \geq \kappa$ and any $\gamma$-sequence $X$ in $\cM$ taking 
  values in $\cB$, the natural map
  $\colim_{\lambda < \gamma} \cM(P,X_{\lambda}) \ra 
                                \cM(P,\colim_{\lambda < \gamma} X_{\lambda})$
  is an isomorphism, where the last colimit is taken in $\cM$.
\end{defn}

\begin{defn}
  Let $I$ be a class of maps in a cocomplete category.
  A map is said to be \mdfn{$I$-injective} if it has the right lifting
  property with respect to each map in $I$, and we write $\Iinj$ for the
  category containing these maps.
  A map is said to be an \mdfn{$I$-cofibration} if it has the left lifting
  property with respect to each map in $\Iinj$, and we write $\Icof$ for the
  category containing these maps.
  A map is said to be \mdfn{$I$-cellular} if it is a transfinite composite
  of pushouts of coproducts of maps in $I$, and we write $\Icell$ for the
  category containing these maps.
\end{defn}

Note that $\Icell$ is a subcategory of $\Icof$.

\begin{defn}
  A \dfn{cofibrantly generated model category} is a model category
  $\cM$ for which there exist sets $I$ and $J$ of morphisms with
  domains which are small relative to $\Icof$ and $\Jcof$, respectively,
  such that $\Icof$ is the category of cofibrations and
  $\Jcof$ is the category of acyclic cofibrations.
  It follows that $\Iinj$ is the category of acyclic fibrations
  and that $\Jinj$ is the category of fibrations.
\end{defn}

For example, take $\cM$ to be the category of spaces and take
$I = \{S^{n}\ra B^{n+1}\}$ and $J = \{B^{n} \times 0 \ra B^{n} \times [0,1]\}$.

\begin{prop} (Recognition Lemma.)\label{pr:rl}
  Let $\cM$ be a category which is complete and cocomplete, let
  $W$ be a class of maps which is closed under retracts and
  satisfies the two-out-of-three axiom, and let $I$ and $J$ be
  sets of maps with domains which are small relative to 
  $\Icell$ and $\Jcell$, respectively, such that
  \begin{roenumerate}
  \item $\Jcell \subseteq \Icof \cap W$ and $\Iinj \subseteq \Jinj \cap W$, and
  \item $\Jcof \supseteq \Icell \cap W$ or $\Iinj \supseteq \Jinj \cap W$.
  \end{roenumerate}
  Then $\cM$ is cofibrantly generated by $I$ and $J$, and $W$ is the
  subcategory of weak equivalences. 
  Moreover, a map is in $\Icof$ if and only if it is a retract
  of a map in $\Icell$, and similarly for $\Jcof$ and $\Jcell$.
  \qed
\end{prop}

The proof, which uses the small object argument, can be found 
in~\cite{dwhika} and~\cite{hi:lmc}.

\section{Derived categories}\label{se:dc}

Let $\A$ be a complete and cocomplete abelian category
equipped with a projective class $\Proj$.
Our goal is to construct a model category structure on the category of
chain complexes of objects of $\A$ in a way that takes into account 
the projective class $\Proj$.
At a certain point we will need a set-theoretic assumption, 
but we proceed as far as possible without it.

We write $\Ch$ for the category of unbounded chain complexes
of objects of $\cA$ and degree zero chain maps.
To fix notation, assume that the differentials lower degree.
For an object $X$ of $\Ch$, define $Z_n X := \ker(d : X_n \ra X_{n-1})$ and
$B_n X := \im(d : X_{n+1} \ra X_n)$, and write $H_n X$ for the quotient.
Given an object $P$ of $\A$ and a chain complex $X$, we write
$\A(P,X)$ for the chain complex which has the abelian group $\A(P,X_{n})$
in degree $n$.
This is the chain complex of $P$-elements of $X$.
(See Section~\ref{se:pc} for a discussion of $P$-elements.)
The cycles, boundaries and homology classes in $\A(P,X)$ are called
\mdfn{$P$-cycles}, \mdfn{$P$-boundaries} and \mdfn{$P$-homology classes} 
in $X$.
Note that a $P$-cycle is the same as a $P$-element of $Z_{n} X$, but
that not every $P$-element of $B_{n} X$ is necessarily a $P$-boundary
(unless $P$ is a categorical projective).

Using this terminology, we can say which maps in our category
should be weak equivalences, fibrations and cofibrations.

\begin{defn}\label{de:ch}
Consider a chain map $f : X \ra Y$.
The map $f$ is said to be a \dfn{quasi-isomorphism} if
it induces an isomorphism $H_{n}X \ra H_{n}Y$ for each $n$.
The map $f$ is said to be a \dfn{weak equivalence}
if for each $\Proj$-projective $P$ the chain map $\A(P,X) \ra \A(P,Y)$
is a quasi-isomorphism, \ie 
if $f$ induces an isomorphism of $P$-homology groups.
The map $f$ is said to be a \dfn{fibration} 
if for each $\Proj$-projective $P$ the chain map $\A(P,X) \ra \A(P,Y)$
is degreewise epic, \ie
if $f$ induces a surjection of $P$-elements.
And $f$ is said to be a \dfn{cofibration} if each $f_{n} : X_{n} \ra Y_{n}$ 
is split-monic and the complex of cokernels is ``cofibrant'': 
a complex is \dfn{cofibrant} if for some ordinal $\kappa$ it is a retract 
of a complex which is a colimit of a $\kappa$-sequence
\[ 0 = C^{0} \lra C^{1} \lra C^{2} \lra \cdots  \]
such that each map is degreewise split-monic
and the successive quotients are complexes of $\Proj$-projective objects 
with zero differential.
As usual, a map is said to be an \dfn{acyclic fibration} 
(resp.\ \dfn{acyclic cofibration})
if it is both a fibration (resp.\ cofibration) and a weak equivalence.
\end{defn}

It will be useful to have the following notation.  
We write $\Sigma$ for the usual automorphism of $\Ch$, which
sends a complex $X$ to the complex $\Sigma X$ with 
$(\Sigma X)_{n} = X_{n-1}$ and $d_{\Sigma X} = -d_X$.
The functor $\Sigma$ is defined on morphisms by $(\Sigma f)_n = f_{n-1}$.
Given an object $A$ in $\cA$, 
we consider $A$ as a chain complex concentrated in degree $0$.  
Thus we write $\Sigma^{n} A$ for the complex with $A$
concentrated in degree $n$.
The maps from $\Sigma^{n} A$ to a complex $X$
biject with the $A$-elements of $Z_{n}X$. 
Similarly, we write $D^{n} A$ for the complex which has $A$ in degrees
$n-1$ and $n$ connected by the identity map, and which is zero in
other degrees,
and we find that 
the maps from $D^{n}A$ to a complex $X$
biject with the $A$-elements of $X_{n}$.
There is a natural map $\Sigma^{n-1} A \ra D^{n} A$.

The following lemma will be used to prove our main result.

\begin{lemma}\label{le:injs}
Consider an object $A \in \A$ and a map $f : X \ra Y$ of chain complexes.
\begin{roenumerate}
\item The map $f$ has the RLP with respect to each of the
maps $0 \ra D^{n}A$, $n \in \Z$, if and only if 
the induced map of $A$-elements is a surjection.
\item The map $f$ has the RLP with respect to each of the
maps $\Sigma^{n-1} A \ra D^{n}A$, $n \in \Z$, if and only if 
the induced map of $A$-elements is a surjection and a quasi-isomorphism.
\end{roenumerate}
\end{lemma}

\begin{proof}
We begin with (\emph{i}): The map $f$ has the RLP with respect to the
map $0 \ra D^{n}A$ if and only if each map $D^{n}A \ra Y$ factors
through $f$, \ie if and only if each $A$-element of $Y_{n}$
is in the image of $f_{n}$.


Now (\emph{ii}): The map $f$ has the RLP with respect
to the map $\Sigma^{n-1}A \ra D^{n}A$ if and only if 
for each $A$-element $y$ of $Y_{n}$ whose boundary is the image of
an $A$-cycle $x$ of $X_{n-1}$, there is an $A$-element $x'$ of 
$X_{n}$ which hits $y$ under $f$ and $x$ under the differential.
(In other words, if and only if the natural map 
$X_{n} \ra Z_{n-1}X \times_{Z_{n-1}Y} Y_{n}$
induces a surjection of $A$-elements.)

So suppose that $f$ has the RLP with respect to each map
$\Sigma^{n-1}A \ra D^{n}A$.
As a preliminary result, 
we prove that $f$ induces a surjection of $A$-cycles.
Suppose we are given an $A$-cycle $y$ of $Y$.
Its boundary is zero and is thus in the image of the $A$-cycle $0$ of $X$.
Therefore $y$ is the image of an $A$-cycle $x'$.

It follows immediately that $f$ induces a surjection in $A$-homology.

Now we prove that $f$ induces a surjection of $A$-elements.
Suppose we are given an $A$-element $y$ of $Y$.
By the above argument, its boundary, which is an $A$-cycle,
is the image of an $A$-cycle $x$ of $X$.
Thus, by the characterization of maps $f$ having the RLP,
we see that there is an $A$-element $x'$ which hits $y$.

Finally, we prove that $f$ induces an injection in $A$-homology.
Suppose $x$ is an $A$-cycle of $X$ whose image under $f$ is
the boundary of an $A$-element $y$.
Then, using the characterization of the RLP, we see that there
is an $A$-element $x'$ whose boundary is $x$.
This shows that $f$ induces an injection in $A$-homology.

We have proved that if $f$ has the RLP with respect to the maps
$\Sigma^{n-1}A \ra D^{n}A$, then $f$ induces an isomorphism in $A$-homology
and a surjection of $A$-elements.
The proof of the converse is similar.
\end{proof}

\begin{cor}\label{co:char}
A map $f : X \ra Y$ is a fibration \ulp resp.\ acyclic fibration\urp\ 
if and only if it has the RLP with respect to the map 
$0 \ra D^{n}P$ \ulp resp.\ $\Sigma^{n-1} P \ra D^{n} P$\urp\ 
for each $\Proj$-projective $P$ and each $n \in \Z$.    \qed
\end{cor}

We want to claim that the above definitions lead to a model
category structure on the category $\Ch$.  In order to
prove this, we need to assume that there is a \emph{set} $\cS$
of $\Proj$-projectives such that a map $f : A \ra B$
is $\Proj$-epic if and only if 
$f$ induces a surjection of $P$-elements for each $P$ in $\cS$.
This implies that for any $B$ there is a $\Proj$-epimorphism $P \ra B$
with $P$ a coproduct of objects from $\cS$, and that every 
$\Proj$-projective object is a retract of such a coproduct.
We also need to assume that each $P$ in $\cS$ is small relative
to the subcategory 
$K := \{(1,0) : A \ra A \oplus Q \st A \in \A \text{ and } Q \in \Proj \}$.
(See the previous section for the definition of ``small.'')
When these conditions hold we say that $\Proj$ is 
\dfn{determined by a set of small objects}.

\begin{thm}\label{th:dc}
Assume that $\Proj$ is determined by a set $\cS$ of small objects.
Then the category $\Ch$ is a cofibrantly generated model category
with the following generating sets:
\[ I := \{ \Sigma^{n-1} P \ra D^{n} P \st P \in \cS,\, n \in \Z \} , \]
\[ J := \{              0 \ra D^{n} P \st P \in \cS,\, n \in \Z \} . \]
Moreover, the weak equivalences, fibrations, cofibrations and cofibrant
objects are as described in Definition~\ref{de:ch}, and every object is
fibrant.
\end{thm}

\begin{proof}
In this proof we use the terms weak equivalence, (acyclic) fibration,
(acyclic) cofibration, and cofibrant as defined in Definition~\ref{de:ch}.

We check the hypotheses of the Recognition Lemma from the previous section.
Since $\cA$ is complete and cocomplete, so is $\Ch$;  limits and
colimits are taken degreewise.
The class of weak equivalences is easily seen to be closed under retracts
and to satisfy the two-out-of-three condition.
The zero chain complex is certainly small relative to $\Jcell$.
We assumed that each $P$ in $\cS$ is small relative to $K \subseteq \cA$, 
and it follows that each $\Sigma^{k} P$ is small relative to 
$\Icell \subseteq \Ch$.
In more detail:  A map is in $\Icell$ if and only if it 
is a degreewise split monomorphism
whose cokernel is a transfinite colimit of degreewise split monomorphisms
whose cokernels are complexes of $\Proj$-projectives with zero differential.
In particular, every map in $\Icell$ is a degreewise split monomorphism
whose cokernel is a complex of $\Proj$-projectives.
Let $P$ be an object of $\cS$.
For large enough regular cardinals $\kappa$,
$P$ is small with respect to $\kappa$-sequences of maps in $K$.
Suppose that for such a $\kappa$ we have a $\kappa$-sequence 
$X^{0} \ra X^{1} \ra \cdots$ of maps in $\Icell$ with colimit $X$.
A chain map from $\Sigma^{k} P$ to $X$ is just a map $f : P \ra X_{k}$
such that $df = 0$.  Given such an $f$, it factors through some
$g : P \ra X^{\alpha}_{k}$, since $P$ is small relative to $K$, and the
sequence $X^{0}_{k} \ra X^{1}_{k} \ra \cdots$ is in $K$.
Moreover, since $df = 0$, it follows that $dg = 0$, at least after
passing to $X^{\beta}_{k}$ for some $\beta \geq \alpha$, again using smallness.
Thus the map $\colim \Ch(\Sigma^{k} P,X^{\lambda}) \ra \Ch(\Sigma^{k} P, X)$
is surjective.
An easier argument shows that it is monic, and therefore $\Sigma^{k} P$
is small relative to $\Icell$.

Since the projectives in $\cS$ may be used to test whether
a map is a fibration or weak equivalence, Corollary~\ref{co:char}
tells us that $\Iinj$ is the collection of acyclic fibrations
and that $\Jinj$ is the collection of fibrations.
Thus we have an equality $\Iinj = \Jinj \cap W$, giving us two
of the inclusions required by the Recognition Lemma.

We now prove that $\Jcell \subseteq \Icof \cap W$.
Since $\Iinj \subseteq \Jinj$, it is clear that $\Jcof \subseteq \Icof$,
so in particular $\Jcell \subseteq \Icof$.
We must prove that $\Jcell \subseteq W$, \ie that each map which
is a transfinite composite of pushouts of coproducts of
maps in $J$ is a weak equivalence.
A map in $J$ is of the form $0 \ra D^{n} P$ for some $P \in \cS$.
Thus a coproduct of maps in $J$ is of the form $0 \ra C$ for
some contractible complex $C$.
A pushout of such a coproduct is of the form $X \ra X \oplus C$,
again with $C$ a contractible complex.
And a transfinite composite of such maps is of the same form as 
well.
In particular, a transfinite composite of pushouts of coproducts
of maps from $J$ is a weak equivalence.

We can now apply the Recognition Lemma and conclude that $\Ch$
is a model category with weak equivalences as described in
Definition~\ref{de:ch}, cofibrations given by the maps in
$\Icof$ and fibrations given by the maps in $\Jinj$.
We saw above that the maps in $\Jinj$ are precisely those
we called fibrations in Definition~\ref{de:ch}, and it is
clear from this that every object is fibrant.

It remains to check that the maps in $\Icof$
are as described in Definition~\ref{de:ch}.
We continue use the words ``cofibration''
and ``cofibrant'' as defined in Definition~\ref{de:ch}.
Suppose that $i$ is in $\Icof$.  By~\ref{pr:rl},
$i$ is a retract of a map $j$ which is a transfinite
composite of pushouts of coproducts of maps in $I$.
As above, one can show that $j$ is degreewise split-monic
and that the complex $C$ of cokernels is cofibrant.
It follows that $i$ is degreewise split-monic and that
the complex of cokernels is a retract of $C$ and is therefore
cofibrant.  That is, $i$ is a cofibration.
And now it is clear that if the map $0 \ra X$ is in $\Icof$,
then $X$ is cofibrant.

Using that $\Icof$ is closed under transfinite compositions and
retracts, one can show that every cofibration is in $\Icof$
and in particular that if $X$ is cofibrant then the map $0 \ra X$ 
is in $\Icof$.
\end{proof}

For the rest of this section, assume that the projective class
$\Proj$ is determined by a set of small objects, so that
$\Ch$ is a model category.
The homotopy category of $\Ch$, formed by inverting the
weak equivalences, is called the \dfn{derived category} of $\A$
(with respect to $\Proj$).  It is denoted $\D$, and a fundamental
result in model category theory asserts that $\D(X,Y)$ is a set
for each $X$ and $Y$.
In Exercise 10.4.5 of \cite{we:iha}, Weibel outlines an argument
which proves that $\D(X,Y)$ is a set when there are enough
(categorical) projectives, $\Proj$ is the categorical projective
class, and $\A$ satisfies AB5.  A connection between Weibel's hypotheses
and ours is that if $\A$ has enough projectives which are
small with respect to all filtered diagrams in $\A$, then AB5 holds.
Of course, the smallness condition needed for our theorem is
weaker than this in that smallness is only required for diagrams
indexed by large cardinals and taking values in the subcategory $K$.

We now show that the notion of homotopy determined from the
model category structure (see \cite{dwsp:htmc}
or \cite{qu:ha}) corresponds to the usual notion of chain homotopy.

\begin{defn}
If $M$ is an object in a model category $\cC$, a \dfn{good cylinder object}
for $M$ is an object $M \times I$ and a factorization
$M \coprod M \llra{i} M \times I \llra{p} M$ of the codiagonal map,
with $i$ a cofibration and $p$ a weak equivalence.
(Despite the notation, $M \times I$ is \emph{not} in general a product 
of $M$ with an object $I$.)
A \dfn{left homotopy} between maps $f, g : M \ra N$ is a map
$H : M \times I \ra N$ such that the composite $H i$ is equal
to $f \coprod g : M \coprod M \ra N$, for some good cylinder
object $M \times I$.
\end{defn}

The notion of \dfn{good path object} $N^{I}$ for $N$ is dual to that of 
good cylinder object and leads to the notion of \dfn{right homotopy}.  
The following standard result can be found in \cite[Section 4]{dwsp:htmc}, 
for example.

\begin{lemma}\label{le:homotopy}
For $M$ cofibrant and $N$ fibrant, two maps $f, g : M \ra N$ 
are left homotopic if and only if they are right homotopic,
and both of these relations are equivalence relations and
respect composition.
Moreover, if $M \times I$ is a fixed good cylinder object
for $M$, then $f$ and $g$ are left homotopic if and only if
they are left homotopic using $M \times I$;  similarly for
a fixed good path object.            \qed
\end{lemma}

Because of the lemma, for $M$ cofibrant and $N$ fibrant we have
a well-defined relation of \dfn{homotopy} on maps $M \ra N$.
One can show that the homotopy category of $\cC$, which is by
definition the category of fractions formed by inverting the
weak equivalences, is equivalent to the category consisting
of objects which are both fibrant and cofibrant with morphisms
being homotopy classes of morphisms.

Now we return to the study of the model category $\Ch$.

\begin{lemma}\label{le:rh}
Let $M$ and $N$ be objects of $\Ch$ with $M$ cofibrant.
Two maps $M \ra N$ are homotopic if and only if they are chain homotopic.
\end{lemma}

\begin{proof}
We construct a factorization $N \ra N^{I} \ra N \times N$ of the
diagonal map $\Delta : N \ra N \times N$ in the following way.
Let $N^{I}$ be the chain complex which has $N_{n} \oplus N_{n+1} \oplus N_{n}$
in degree $n$.  We describe the differential by saying that it sends
a generalized element $(n,\bar{n},n')$ in $(N^{I})_{n}$ to 
$(dn,n-n'-d \bar{n},dn')$.
Let $\alpha : N \ra N^{I}$ be the map which sends $n$ to $(n,0,n)$
and let $\beta : N^{I} \ra N \times N$ be the map which sends
$(n,\bar{n},n')$ to $(n,n')$.
One can check easily that $N^{I}$ is a chain complex and that
$\alpha$ and $\beta$ are chain maps whose composite is $\Delta$.
The map $\alpha$ is a chain homotopy equivalence with chain homotopy inverse
sending $(n,\bar{n},n')$ to $n$; this implies
that it induces a chain homotopy equivalence of generalized elements
and is thus a weak equivalence.
The map $\beta$ is degreewise split-epi; this implies that 
it induces a split epimorphism of generalized elements and 
is thus a fibration.
Therefore $N^{I}$ is a good path object for $N$.  

It is easy to see that a chain homotopy between two maps $M \ra N$
is the same as a right homotopy using the path object $N^{I}$.
By Lemma~\ref{le:homotopy}, 
two maps are homotopic if and only if they are right homotopic using $N^{I}$. 
Thus the model category notion of homotopy is the same
as the notion of chain homotopy when the source is cofibrant.
\end{proof}

There is a dual proof which proceeds by constructing a specific good cylinder
object $M \times I$ for $M$ such that a left homotopy using $M \times I$
is the same as a chain homotopy.

\begin{cor}\label{co:ext}
Let $A$ and $B$ be objects of $\A$ considered as chain complexes
concentrated in degree $0$.
Then $\D(A,\Sigma^{n} B) \iso \Ext_{\Proj}^{n}(A,B)$.
\end{cor}

Here we denote by $\Ext_{\Proj}^{n}(-,B)$ the right derived
functors of $\A(-,B)$ with respect to the projective class $\Proj$.
The group $\Ext_{\Proj}^{n}(A,B)$ can be calculated by forming
a $\Proj$-projective $\Proj$-exact resolution of $A$, applying
$\A(-,B)$, and taking cohomology.
One can show that $\Ext_{\Proj}^{n}(A,B)$ classifies
equivalence classes of $\Proj$-exact sequences
$0 \ra B \ra C_{1} \ra \cdots \ra C_{n} \ra A \ra 0$.

\begin{proof}
The group $\D(A,\Sigma^{n}B)$ may be calculated by choosing a
cofibrant replacement $A'$ for $A$ and computing the homotopy
classes of maps from $A'$ to $\Sigma^{n} B$.  
(Recall that all objects are fibrant, so there is no need to take a 
fibrant replacement for $\Sigma^{n} B$.)
A $\Proj$-projective $\Proj$-exact resolution $P$ of $A$ serves as a 
cofibrant replacement for $A$, and by Lemma~\ref{le:rh} the homotopy 
relation on $\Ch(P,\Sigma^{n}B)$ is chain homotopy, so it follows that
$\D(A,\Sigma^{n}B)$ is isomorphic to $\Ext_{\Proj}^{n}(A,B)$.
\end{proof}

More generally, a similar argument shows that the derived functors
of a functor $F$ can be expressed as the cohomology of the derived
functor of $F$ in the model category sense.  To make the story 
complete, we next show that the shift functor $\Sigma$ corresponds
to the notion of suspension that the category $\D$ obtains as the
homotopy category of a pointed model category.

\begin{defn}\label{de:susp}
Let $\cC$ be a pointed model category.
If $M$ is cofibrant, we define the \dfn{suspension} $\Sigma M$ of $M$ to
be the cofibre of the map $M \coprod M \ra M \times I$ for some
good cylinder object $M \times I$.
(The cofibre of a map $X \ra Y$ is the pushout $* \coprod_X Y$, 
where $*$ is the zero object.)
$\Sigma M$ is cofibrant and well-defined up to homotopy equivalence.
\end{defn}

The loop object $\Omega N$ of a fibrant object $N$ is defined dually.
These operations induce adjoint functors on the homotopy category.
A straightforward argument based on the path object described above
(and a dual cylinder object) proves the following lemma.

\begin{lemma}\label{le:loop}
In the model category $\Ch$,
the functor $\Sigma$ defined in Definition~\ref{de:susp} can be
taken to be the usual suspension, so that
$(\Sigma X)_{n} = X_{n-1}$ and $d_{\Sigma X} = -d_{X}$.
Similarly, $\Omega X$ can be taken to be
the complex $\Sigma^{-1}X$.  That is, $(\Omega X)_{n} = X_{n+1}$
and $d_{\Omega X} = -d_{X}$.  \qed
\end{lemma}

In particular, $\Sigma$ and $\Omega$ are inverse functors.
Hovey~\cite{ho:mc} has shown that this implies that cofibre
sequences and fibre sequences agree (up to the usual sign)
and that $\Sigma$ and the cofibre sequences give rise to a 
triangulation of the homotopy category.
(See~\cite[Section I.3]{qu:ha} for the definition of cofibre and
fibre sequences in any pointed model category.)
Applying Hovey's result to our situation we can be more explicit.
Given a map $f : L \ra M$ of chain complexes, the \dfn{standard triangle}
on $f$ is the sequence $L \ra M \ra N \ra \Sigma L$, where 
$N_{n} = M_{n} \oplus L_{n-1}$, $d_{N}(m,l) = (d_{M} m + f l, -d_{L} l)$ on 
generalized elements, and the maps are the inclusion and projection.

\begin{cor}\label{co:triangulated}
The category $\D$ is triangulated.
A sequence $L' \ra M' \ra N' \ra \Sigma L'$ is a triangle if and only if
it is isomorphic in $\D$ to a standard triangle.
\end{cor}

\begin{proof}
That $\D$ is triangulated follows from Hovey's result, since we have shown
that $\Sigma$ is a self-equivalence.
The identification of the triangles uses the definition of fibre sequences 
from~\cite{qu:ha} and the construction of the path object from the proof 
of Lemma~\ref{le:rh}.
There is also a dual proof using cofibre sequences and cylinder objects.
\end{proof}

The results of this section probably apply with $\cA$ replaced by
a complete and cocomplete additive category.

\section{The pure and categorical derived categories}\label{se:examples}

In this section we let $R$ be an associative ring with unit and 
we take for $\A$ the category of left $R$-modules.
We are concerned with two projective classes on the category $\A$.
The first is the categorical projective class $\cC$ 
whose projectives are summands of free modules, 
whose exact sequences are the usual exact sequences,
and whose epimorphisms are the surjections.
The second is the pure projective class $\cP$ 
whose projectives are summands of sums of finitely presented modules.
A short exact sequence is $\Proj$-exact iff it remains exact
after tensoring with any right module.
A map is a $\Proj$-epimorphism iff it appears in a $\Proj$-exact
short exact sequence.
We say pure projective instead of $\Proj$-projective, and similarly
for pure exact and pure epimorphism.
We assume that the reader has some familiarity with these
projective classes.  A brief summary with further references
may be found in~\cite[Section~9]{ch:itcpgs}.  As usual, we write
$\Ext^{*}(-,B)$ (resp.\ $\PExt^{*}(-,B)$) for the derived functors of 
$\A(-,B)$ with respect to the categorical (resp.\ pure) projective 
class. 

Both of these projective classes are determined by sets of small
objects:  $\cC$ is determined by $\{R\}$ and $\cP$ is determined
by any set of finitely presented modules containing a representative
from each isomorphism class.
Thus we get two model category structures on $\Ch$ and two
derived categories, the categorical derived category $\DC$ and
the pure derived category $\DP$.
We refer to the pure weak equivalences in $\Ch$ as pure quasi-isomorphisms,
and as usual call the categorical weak equivalences simply quasi-isomorphisms.
Similarly, we talk of pure fibrations and fibrations, pure cofibrations
and cofibrations, etc.

Pure homological algebra is of interest in stable homotopy
theory because of the following result.

\begin{thm}\cite{chst:pmht}
Phantom maps from a spectrum $X$ to an Eilenberg-Mac\,Lane spectrum $HG$
are given by $\PExt_{Z}^{1}(H_{-1} X, G)$, that is, by maps of degree
one from $H_{-1} X$ to $G$ in the pure derived category of abelian groups.
\end{thm}

The two derived categories are connected in various ways.
For example, since every categorical projective is pure projective, it follows
that every pure quasi-iso\-mor\-phism is a quasi-isomorphism.  This
implies that there is a unique functor $R: \DP \ra \DC$
commuting with the functors from $\Ch$.  This functor has a left
adjoint $L$.  We can see this in the following way:  the identity
functor on $\Ch$ is adjoint to itself.  It is easy to see that
every pure fibration is a fibration and that every pure acyclic
fibration is an acyclic fibration.  Therefore, 
by~\cite[Theorem~9.7]{dwsp:htmc}, there is an induced pair
of adjoint functors between the pure derived category and
the categorical derived category, and the right adjoint is the functor $R$
mentioned above. 
This right adjoint is the identity on objects, since everything is
fibrant, and induces the natural map $\PExt^{*}(A,B) \ra \Ext^{*}(A,B)$
for $R$-modules $A$ and $B$.
The left adjoint sends a complex $X$ in $\DC$ to
a (categorical) cofibrant replacement $\tilde{X}$ for $X$.
Similar adjoint functors exist whenever one has two projective
classes on a category, one containing the other.

In addition to the connection between phantom maps and pure homological
algebra, the author is interested in the pure derived category as a
tool for connecting the global pure dimension of a ring $R$ to
the behaviour of phantom maps in $\DC$ and $\DP$ under composition.
This work will hopefully appear in a future paper.

\section{Simplicial objects and the bounded below derived category}\label{se:simp}

In this section, let $\cA$ be a complete and cocomplete pointed category
with a projective class $\Proj$.  We do not assume that $\cA$ is abelian.
Write $\sA$ for the category of simplicial objects in $\cA$.
Given a simplicial object $X$ and an object $P$, write $\cA(P,X)$
for the simplicial set which has $\cA(P,X_{n})$ in degree $n$.
Consider the following conditions:
\begin{itemize}
\item [(*)] For each $X$ in $\sA$ and each $P$ in $\Proj$, $\cA(P,X)$
is a fibrant simplicial set.
\item [(**)] $\Proj$ is determined by a set $\cS$ of small objects.
Here we require that each $P$ in $\cS$ be small with respect to all
monomorphisms in $\cA$, not just the split monomorphisms with
$\Proj$-projective cokernels.
\end{itemize}

\begin{thm}\label{th:simp}
Let $\cA$ be a complete and cocomplete pointed category with a 
projective class $\Proj$.  If $\cA$ and $\Proj$ satisfy \ulp *\urp\ or 
\ulp **\urp, then
$\sA$ is a simplicial model category.  The weak equivalences
\ulp resp.\ fibrations\urp\ are the maps $f$ such that $\cA(P,f)$
is a weak equivalence \ulp resp.\ fibration\urp\ of simplicial sets
for each $P$ in $\cS$.
(Here we write $\cS$ for the set of small objects which determines
$\Proj$ when \ulp **\urp\ holds, and we take $\cS$ to be equal to
$\Proj$ when \ulp *\urp\ holds.)
In general the cofibrations are determined by the lifting property.
When \ulp **\urp\ holds, $\sA$ is cofibrantly generated by the sets
\[ I := \{ P \tensor \dot{\Delta}[n] \ra P \tensor \Delta[n] 
       \st P \in \cS,\, n \geq 0 \} , \]
\[ J := \{ P \tensor V[n,k] \ra P \tensor \Delta[n] 
       \st P \in \cS,\, 0 < n \geq k \geq 0 \} . \]
\end{thm}

In the above, $\Delta[n]$ denotes the standard $n$-simplex simplicial set;  
$\dot{\Delta}[n]$ denotes its boundary;
and $V[n,k]$ denotes the subcomplex with the $n$-cell and its
$k$th face removed.
For an object $P$ of $\cA$ and a set $K$, $P \tensor K$ denotes
the coproduct of copies of $P$ indexed by $K$.
When $K$ is a simplicial set, $P \tensor K$ denotes the simplicial
object in $\cA$ which has $P \tensor K_{n}$ in degree $n$.

If $\cA$ is abelian, then (*) holds.  Moreover, in this case
there is an equivalence of categories given by
the normalization functor $\sA \ra \Ch^{+}$,
where $\Ch^{+}$ denotes the category of non-negatively graded
chain complexes of objects of $\cA$.  
Thus we can deduce:

\begin{cor}
Let $\cA$ be a complete and cocomplete abelian category with
a projective class $\Proj$.  Then $\Ch^{+}$ is a model category.
A map $f$ is a weak equivalence iff $\cA(P,f)$ is
a quasi-isomorphism for each $P$ in $\Proj$.
A map $f$ is a fibration iff $\cA(P,f)$ is surjective in positive
degrees \ulp but not necessarily in degree 0\urp\ for each $P$ in $\Proj$.
A map is a cofibration iff it is degreewise split-monic with degreewise
$\Proj$-projective cokernel. 
Every complex is fibrant, and a complex is cofibrant iff
it is a complex of $\Proj$-projectives. \qed
\end{cor}

Note that no conditions on the projective class are required, and that
the description of the cofibrations is simpler than in the unbounded case.

As another example of the theorem, let $G$ be a group and consider
$G$-simplicial sets, or equivalently, simplicial objects
in the category $\cA$ of (say, left) $G$-sets.
Let $\cF$ be a family of subgroups of $G$, and consider the
set $\cS = \{ G/H \st H \in \cF \}$ of homogeneous spaces.
These are small, and determine a projective class $\Proj$.
Thus, using case (**) of the above theorem, we can deduce:

\begin{cor}
Let $G$ be a group and let $\cF$ be a family of subgroups of $G$.
The category of $G$-equivariant simplicial sets has a model
category structure in which the weak equivalences are precisely
the maps which induce a weak equivalence on $H$-fixed points
for each $H$ in $\cF$.  \qed
\end{cor}

We omit the proof of Theorem~\ref{th:simp}, and simply note that
it follows the argument in \cite[Section II.4]{qu:ha} fairly
closely.  It is a bit simpler, in that Quillen spends
part of the time (specifically Proposition 2) proving that 
effective epimorphisms give rise to a projective class,
although he doesn't use this terminology.
It is also a bit more complicated, in the (**) case, in that a
transfinite version of Kan's $\Ex^{\infty}$ functor~\cite{ka:exi}
is required, since we make a weaker smallness assumption.
Quillen's argument in this case can be interpreted as a verification
of the hypotheses of the recognition lemma for cofibrantly generated
model categories (Proposition~\ref{pr:rl}).

\bibliography{jourabbrev,christensen}
\bibliographystyle{plain}

\end{document}